\documentstyle[11pt]{article}
\font\tenmath=msbm10 scaled 1200
\font\sevenmath=msbm7 scaled 1200
\font\fivemath=msbm5 scaled 1200

\newfam\mathfam \textfont\mathfam=\tenmath
\scriptfont\mathfam=\sevenmath \scriptscriptfont\mathfam=\fivemath
\def\math{\fam\mathfam}
\def\R{{\math R}}
\def\N{{\math N}}

\def \^#1{\if#1i{\accent"5E\i}\else{\accent"5E#1}\fi}

\def \cqfd{\quad_\diamondsuit}

\newtheorem{Thm}{Theorem}
\newtheorem{Lem}{Lemma}

\newtheorem{Cor}{Corollary}

\author{\sc Gilles Pag\`es
\thanks{Laboratoire de Probabilit\'es et Mod\'elisation al\'eatoire, UMR~7599, Universit\'e Paris 6, case 188, 4,
pl. Jussieu, F-75252 Paris Cedex 5. E-mail: {\tt  gpa@ccr.jussieu.fr}}} 

\title{\bf A two armed bandit type problem  revisited} 

\begin{document}

\maketitle
\begin{abstract}
In~\cite{BENBEN} M. Bena\"{\i}m and G. Ben Arous  solve a multi-armed bandit problem arising in the theory of learning in games. We propose
an short elementary  proof of this result  based on a variant of the Kronecker Lemma. 
\end{abstract}




\bigskip

\bigskip

\bigskip
In~\cite{BENBEN}  a multi-armed bandit problem is addressed and investigated by M. Bena\"{\i}m and G. Ben Arous. Let $f_0,\ldots,f_d$
denote
$d+1$ real-valued continuous functions defined on
$[0,1]^{d+1}$.  Given a   sequence $x = (x_n)_{n\ge 1}\!\in\{0,\ldots,d\}^{\N^*}$ (the {\em strategy}),  set  for every $n\ge
1$
\[
\bar x_n := (\bar x_n^0,\bar x_n^1, \ldots,\bar x_n^d)\quad\mbox{ with }\quad \bar x^i_n := \frac{1}{n}\sum_{k=1}^n {\mathbf
1}_{\{x_k=i\}},\;\; i=0,\ldots,d, 
\]  
and
\[
Q(x) = \displaystyle  \liminf_{n\to +\infty} \frac{1}{n} \sum_{k=0}^{n-1} f_{x_{k+1}}(\bar x_{k}).
\]
($\bar
x_0 := (\bar x_0^0,\bar x_0^1, \ldots,\bar x_0^d) \!\in [0,1]^{d+1}$, $\bar x^0_0+\cdots +\bar x^{d}_0=1$ is a starting distribution). Imagine
$d+1$ players enrolled in a cooperative/competitive game with the following simple rules: if player $i\!\in\{0,\ldots,d\}$ plays at time
$n$ he is rewarded by $f_{i}(\bar x_n)$, otherwise he gets nothing; only one player can play at the same time. Then the sequence  $x$ is a
playing strategy for the group of players and $Q(x)$ is the {\em global} cumulative worst payoff rate of the strategy $x$ for the whole
community of players (regardless of the cumulative payoff rate  of each player). 

In~\cite{BENBEN} an  answer (see~Theorem~\ref{Ben^2} below) is provided to the following question

\medskip
\centerline{\em What are the good strategies (for the group)~?}

\medskip
\noindent The authors rely on some recent tools developed in  stochastic approximation theory (see $e.g.$~\cite{BEN}). The aim of this  note
is  to provide an elementary and shorter proof based on a slight  improvement of the Kronecker Lemma.

\medskip
 Let ${\cal S}_d:=\{v\!\in[0,1]^d,\, \sum_{i=1}^d v_i \le 1\}$ and ${\cal P}_{d+1}:= \{u\!\in[0,1]^{d+1},\, \sum_{i=1}^{d+1} u_i = 1\}$.
Furthermore, for notational convenience,  set
\begin{eqnarray*}
\forall\, v=(v_1,\ldots,v_d)\in {\cal S}_d,\quad \tilde v&:=& (1-\sum_{i=1}^d v_i,
v_1,\ldots,v_d)\in {\cal P}_{d+1},\\
\forall\,u=(u_0,u_1,\ldots,u_d)\!\in {\cal P}_{d+1}, \quad \hat u &:=&(u_1,\ldots, u_d)\in {\cal S}_d. 
\end{eqnarray*} 
The canonical inner product on $\R^d$ will be  denoted by   $(v|v')=\sum_{i=1}^d v_iv'_i$. The interior of a subset $A$ of $\R^d$ will
be denoted $ A^{\hskip -0.15 cm^{^\circ}}$. For a sequence $u=(u_n)_{n\ge 1}$,  $\Delta u_n \!:=\!u_n-u_{n-1}$, $n\!\ge\!1$.

\medskip
The main result   is the following theorem (first established in~\cite{BENBEN}). 
\begin{Thm}\label{Ben^2} Assume   there is  a function $\Phi:{\cal S}_d \to
\R$, continuously differentiable on $ {\cal S }^{\hskip -0.15 cm^{^\circ}}_d$ having  a continuous extension $\nabla \Phi$ on  ${\cal S}_d$ and
satisfying:
\begin{equation}\label{gradPhi}
\forall\, v\in {\cal S}_d,\quad \nabla \Phi(v) = \left(f_i(\tilde v) -f_0(\tilde v) \right)_{1\le i\le d}.
\end{equation}    
Set for every $u\!\in{\cal P}_{d+1}$, 
\[
q(u):= \sum_{i=0}^{d+1} u_i\, f_i(u)
\]
and $Q^* := \max \left\{q(u), \; u\!\in{\cal P}_{d+1}\right\}$. Then, for every strategy $x\!\in
\{0,1,\ldots,d\}^{\N^*}$, 
\[
Q(x)\le Q^*. 
\]
Furthermore, for any strategy $x$ such that $\bar x_n \to \bar x_\infty$, 
\[
\frac{1}{n}\sum_{k=1}^n f_{x_{k+1}}(\bar x_k) \to q(\bar x_\infty)\quad \mbox{ as }\quad n\to \infty \qquad (\hbox{so that } Q(x) = q(\bar
x_\infty)).
\] 
In particular there is no better strategy than choosing the player  at random according to an i.i.d. strategy   with distribution
$\bar x^*\!\in {\rm argmax}\,q$. 
\end{Thm}

The key of the proof is the following slight extension of the Kronecker Lemma. 
\begin{Lem}\label{Un} (``\`a la Kronecker" Lemma) Let $(b_n)_{n\ge 1}$ be a nondecreasing sequence of positive real numbers converging to
$+\infty$ and let
$(a_n)_{n\ge 1}$ be a sequence of real numbers. Then
\[
 \liminf_{n\to +\infty} \sum_{k=1}^n\frac{a_k}{b_k} \in  \R \quad\Longrightarrow\quad   \liminf_{n\to +\infty}\frac{1}{b_n}
\sum_{k=1}^na_k\le 0.
\] 

\end{Lem}

\noindent {\bf Proof.} Set $\displaystyle C_n= \sum_{k=1}^n \frac{a_k}{b_k},\; n\ge 1$ and $C_0=0$ so that $a_n =  b_n \Delta C_n$.
As a consequence, an Abel transform yields
\begin{eqnarray*}
\frac{1}{b_n}Ê\sum_{k=1}^n a_k &=& \frac{1}{b_n}Ê\sum_{k=1}^n b_k\Delta C_k=\frac{1}{b_n}Ê\left(b_nC_n - \sum_{k=1}^n C_{k-1}\Delta
b_k\right)
\\ &=& C_n - \frac{1}{b_n} \sum_{k=1}^n C_{k-1}\Delta b_k. 
\end{eqnarray*}
Now, $\displaystyle  \liminf_{n\to +\infty} C_n$ being finite, for every $\varepsilon>0$, there is an integer $n_\varepsilon$ such that for every $k\ge
n_{\varepsilon}
$,
$ C_k \ge \displaystyle  \liminf_{n\to +\infty} C_n -\varepsilon$. Hence
\[
\frac{1}{b_n} \sum_{k=1}^n C_{k-1}\Delta b_k\ge \frac{1}{b_n} \sum_{k=1}^{n_\varepsilon} C_{k-1}\Delta b_k+
\frac{b_n-b_{n_\varepsilon}}{b_n}\left(\displaystyle  \liminf_kC_k-\varepsilon\right).
\]
Consequently, $\displaystyle\displaystyle  \liminf_{n\to +\infty} C_n$ being finite, one concludes that  
\begin{eqnarray*}
\displaystyle  \liminf_{n\to +\infty} \frac{1}{b_n}Ê\sum_{k=1}^n a_k&\le& \displaystyle  \liminf_{n\to +\infty} C_n -0- 1\times \left(\displaystyle  \liminf_{k\to
+\infty}C_k-\varepsilon\right)= \varepsilon.\qquad \cqfd
\end{eqnarray*}

\noindent {\bf Proof of Theorem~\ref{Ben^2}.} First note that for every $u= (u_0,\ldots,u_d)\!\in{\cal P}_{d+1}$, 
\begin{eqnarray*}
q(u) := \sum_{i=0}^{d+1} u_i f_i(u)= f_0(u) +\sum_{i=1}^d u_i(f_i(u)-f_0(u))
\end{eqnarray*}
so that
\begin{eqnarray*}
Q^* =  \sup_{v\in {\cal S}_d}\left\{f_0(\tilde v)+ \sum_{i=1}^d v_i(f_i(\tilde v)-f_0(\tilde v))\right\}=  \sup_{v\in {\cal
S}_d}\left\{f_0(\tilde v)+
\left(v|\nabla \Phi(v)\right)\right\}.
\end{eqnarray*}
Now, for every $k\ge 0$ 
\begin{eqnarray*}
f_{x_{k+1}}(\bar x_k)-q(\bar x_k)& = &\sum_{i=0}^d(f_i(\bar x_k) {\mathbf1}_{ \{x_{k+1}=i\}}-\bar x_k^i f_i(\bar x_k))= \sum_{i=0}^df_i(\bar
x_k)( {\mathbf1}_{ \{x_{k+1}=i\}}-\bar x_k^i)\\ &=& \sum_{i=0}^df_i(\bar x_k)(k+1) \Delta \bar x^i_{k+1}\\
&=& (k+1)\sum_{i=1}^d(f_i(\bar x_k) -f_0(\bar x_k))\Delta \bar x^i_{k+1}.
\end{eqnarray*}
The last equality   reads using Assumption~(\ref{gradPhi}), 
\begin{eqnarray*} 
f_{x_{k+1}}(\bar x_k)-q(\bar x_k) &=&(k+1) (\nabla \Phi(\hat{\bar x}_k)|\Delta \hat{\bar x}_{k+1})
\end{eqnarray*}
Consequently, by the fundamental formula of calculus applied to $\Phi$ on $( \hat{\bar x}_{k}, \hat{\bar x}_{k+1})\subset  {\cal S }^{\hskip -0.15
cm^{^\circ}}_d$,  
\begin{eqnarray*}
\frac{1}{n}\sum_{k=0}^{n-1}f_{x_{k+1}}(\bar x_k)-q(\bar x_k)&=&\frac{1}{n}\sum_{k=0}^{n-1}(k+1) \left(\Phi(\hat{\bar x}_{k+1})-\Phi(\hat{\bar
x}_k)\right)-R_n\\
\mbox{with }\hskip 2,5 cm  R_n &:=&\frac{1}{n}\sum_{k=0}^{n-1}\left (\nabla \Phi(\hat \xi_k)-\nabla \Phi(\hat{\bar x}_k)|(k+1)\Delta \hat{\bar
x}_{k+1}\right)   \hskip 2 cm 
\end{eqnarray*}
and $\hat \xi_k \!\in(\hat{\bar x}_k,\hat{\bar x}_{k+1}) ,\; k=1,\ldots n$. The fact  that $|(k+1)\Delta \hat{\bar x}_{k+1}|\le 1$   implies 
\[
|R_n|\le \frac{1}{n} \sum_{k=0}^{n-1} w(\nabla  \Phi,|\Delta \hat{\bar x}_{k+1}|) 
\]
where $w(g,\delta)$ denotes the uniform continuity $\delta$-modulus of   a function $g$. One derives from the uniform continuity of $\nabla
\Phi$ on the compact set ${\cal S}_d$   that 
\[
R_n\to 0 \quad \mbox{ as   } \quad n\to +\infty.
\]
Finally, the continuous function $\Phi$ being bounded on the compact set ${\cal S}_d$, the  partial sums 
\[
  \sum_{k=0}^{n-1} \Phi(\hat{\bar x}_{k+1})-\Phi(\hat{\bar x}_k) = \Phi(\hat{\bar x}_{n+1})-\Phi(\hat{\bar x}_0)
\] 
remain bounded as $n$ goes to infinity. Lemma~\ref{Un} then implies that 
\[
\displaystyle  \liminf_{n\to +\infty} \frac{1}{n}\sum_{k=0}^{n-1}(k+1) \left(\Phi(\hat{\bar x}_{k+1})-\Phi(\hat{\bar x}_k)\right)\le 0.
\]
One concludes by noting that on one hand   
\[
\limsup_{n\to \infty} \frac{1}{n}\sum_{k=0}^{n-1} q(\bar x_k)\le Q^*=\sup_{{\cal P}_{d+1}} q
\]
and that, on the other hand, the function $q$ being continuous, 
\[
\lim_{n\to \infty} \frac{1}{n}\sum_{k=0}^{n-1}q(\bar x_k)= q(x^*)\quad \mbox{ as soon as } \quad \bar x_n \to x^*.\cqfd
\]

 \begin{Cor} When  $d+1=2$ (two players), Assumption~(\ref{gradPhi}) is  satisfied as soon as $f_0$ and $f_1$ are continuous on ${\cal
P}_2$ and then the conclusions of  Theorem\ref{Ben^2} hold true.  
\end{Cor}

\noindent {\bf Proof:} This follows from the obvious fact that the continuous function $u_1\mapsto f_1(1-u_1,u_1)-f_0(1-u_1,u_1)$ on $[0,1]$
  has an antiderivative. $\cqfd$


\bigskip
\noindent {\sc Further comments:} $\bullet$ If one considers a slightly more general game in which some {\em weighted strategies} are allowed,
the final result is not modified in any way provided   the weight sequence satisfies a very light assumption. Namely, assume that at time $n$
the reward  is 
\[
\Delta_{n+1} f_{x_{n+1}}(\bar x_n) \hskip 1 cm \mbox{ instead of }\hskip 1 cm f_{x_{n+1}}(\bar x_n)
\] 
where the weight sequence $\Delta =(\Delta _n)_{n\ge 1}$ satisfies
\[
\Delta_n \ge 0,\; n\ge 1,\quad S_n = \sum_{k= 1}^n \Delta_k \to +\infty,\quad \frac{\Delta_n}{S_n} \to 0\mbox{ as } \; n\to \infty 
\]
then the quantities $\bar x^{\Delta}_0\!\in{\cal P}_{d+1}$, 
$\displaystyle  \bar x^{\Delta}_n := (\bar x^{\Delta,0}_n,\ldots,\bar x^{\Delta,d}_n)$ with $\bar x^{\Delta,i}_n =
\frac{1}{S_n}\sum_{k=1}^n
\Delta_k{\mathbf 1}_{\{x_k=i\}} ,\; i=0,\ldots,d,\; n\ge 1$, 
and $\displaystyle 
Q^{\Delta}(x) = \displaystyle  \liminf_{n\to +\infty} \frac{1}{S_{n}} \sum_{k=0}^{n-1} \Delta_{k+1}f_{x_{k+1}}(\bar x^{\Delta}_{k})$
satisfy all the conclusions of Theorem~\ref{Ben^2} {\em mutatis mutandis}.

\bigskip
$\bullet$ Several applications of  Theorem~\ref{Ben^2} to the theory of learning in games and to stochastic fictitious play are
extensively    investigated  in~\cite{BENBEN} which  we refer to for all these aspects.  As far as we are concerned we will simply make   a
remark about some  ``natural"  strategies which illustrates the theorem in an elementary way. 

In the reward function  at time $k$, $i.e.$ $f_{x_k}(\bar x_{k-1})$, $x_{k}$ represents the competitive term (``who will play~?") and
$\bar x_{k-1}$ represents a cooperative term (everybody's past behaviour  has influence on everybody's reward).

This cooperative/competitive antagonism induces that  in such a game a {\em greedy}  competitive strategy is usually  not optimal (when the
players do not play a symmetric r\^ole). Let us be more specific. Assume for the sake of simplicity that
$d+1=2$ (two players). Then one may consider without loss of generality that $\bar x_n = \hat{\bar x}_n$ $i.e.$ that $\bar x_n$ is a
$[0,1]$-valued real number.  A {\em greedy  competitive}  strategy is defined  by 
\begin{equation}\label{greedy}
\hbox{player $1$ plays at time $n\;$  ($i.e.$ $x_n=1$) iff }\; f_1(\bar x_{n-1})\ge f_0(\bar x_{n-1}) 
\end{equation}  
$i.e.$ the player with the highest reward is nominated to play. Note that such a strategy  is    anticipative from a
probabilistic viewpoint. Then, for every $n\ge 1$, 
\[
f_{x_n}(\bar x_{n-1}) = \max(f_0(\bar x_{n-1}),f_1(\bar x_{n-1})) 
\]
and it is clear that
\[
 f_{x_n}(\bar x_{n-1})- q(\bar x_n) =\max(f_0(\bar x_{n-1}),f_1(\bar x_{n-1}))- q(\bar x_n) =:\varphi(\bar x_n)\ge 0.
\]
On the other hand, the proof of Theorem~\ref{Ben^2} implies that 
\[
  \liminf_{n\to +\infty} \frac{1}{n} \sum_{k=0}^{n-1} \varphi(\bar x_n) \le 0.
\]
Hence, there is at least one weak limiting distribution $\bar \mu_\infty$ of the sequence of empirical measures $\bar \mu_n:= \frac 1n
\sum_{0\le k\le n-1}
\delta_{\bar x_k}$ which is supported by   the closed set
$\{\varphi =0\}\subset \{0, 1\}\cup \{f_0=f_1\}$; on the other ${\rm supp} (\mu_\infty)$ is contained in the set  $\bar {\cal X}_\infty$
of  the limiting values of the sequence $(\bar x_n)$ itself (in fact $\bar {\cal X}_\infty$ is an interval since $(\bar x_n)_n$ is
bounded and
$\bar x_{n+1}-
\bar x_n\to 0$). Hence $\bar {\cal X}_\infty\cap (\{0, 1\}\cup \{f_0=f_1\}) \neq \emptyset$. 

If the greedy strategy $(\bar x_n)_n$ is optimal then ${\rm dist}(\bar x_n, {\rm argmax} \, q)\to 0$ as $n\to \infty$ $i.e.$ $\bar {\cal
X}_\infty \subset {\rm argmax} \, q$. Consequently if
\begin{equation}
{\rm argmax} \, q \cap (\{0,1\}\cup\{f_0=f_1\}) = \emptyset
\end{equation}
then {\em the purely competitive strategy is never   optimal}. 

\medskip
So is the case if 
\[
f_0(x) = a\,x\qquad \mbox{ and }\qquad f_1(x) =b\,(1-x),\qquad x\!\in[0,1],
\]
for some positive parameters $a\neq b$, then 
\[
{\rm argmax} q = \{1/2\} \qquad \mbox{ and }\qquad  f_0(1/2)\neq f_1(1/2).
\] 

In fact, one shows that the greedy strategy $x=(x_n)_{n\ge 1}$ defined by~(\ref{greedy}) satisfies
\[
\bar x_n \to \frac{b}{a+b} \qquad \mbox{and} \qquad   Q(x)=\frac{ab}{a+b}\quad\mbox{ as }\quad n\to \infty
\]
whereas any optimal (cooperative) strategy (like the $i.i.d.$ Bernoulli($1/2$) one) yields an asymptotic (relative) global payoff rate
\[
Q^*= \max_{[0,1]} q = \frac{a+b}{4}.
\]
Note that $Q^*>\frac{ab}{a+b}$ since $ a\neq b$.  (When $a=b$ the greedy strategy becomes optimal.)

\bigskip
$\bullet$  A more  abstract version of Theorem~\ref{Ben^2} can be established using the same approach. The finite set
$\{0,1,\ldots,d\}$ is replaced by a compact metric set $K$, ${\cal P}_{d+1}$ is replaced by the convex set ${\cal P}_{_{\!K}}$ of  
probability distributions on $K$ equipped with the weak topology and the  continuous function $f:K\times {\cal P}_K\to \R$ still derives from a
potential function in some sense. 


\begin{thebibliography}{}
\bibitem{BEN} {\sc M. Bena\"{\i}m} (1999). Dynamics of stochastic algorithms, in J. Az\'ema et al. eds, {\em S\'eminaire de
probabilit\'es XXXIII}, L.N. in Math. 1708, 1-68, Springer  Verlag, Berlin.
\bibitem{BENBEN} {\sc M. Bena\"{\i}m, G. Ben Arous} (2003). A two armed
bandit type problem, {\em Game Theory}, {\bf 32}(3), 3-16. 

\end{thebibliography}
\end{document}